\documentclass[11pt]{article}
\usepackage{}
\usepackage{mathrsfs}
\usepackage{amsfonts}

\usepackage{CJK}
\usepackage{amsfonts,amssymb,amsmath,mathrsfs}

\usepackage{color,latexsym,amsfonts}
 \newtheorem{theorem}{Theorem}[section]
 \newtheorem{lemma}[theorem]{Lemma}
 \newtheorem{corollary}[theorem]{Corollary}
 \newtheorem{proposition}[theorem]{Proposition}
 
 \newtheorem{Definition}[theorem]{Definition}
 \newtheorem{remark}[theorem]{Remark}
 \newtheorem{condition}[theorem]{Condition}

 \def\blemma{\begin{lemma}\sl{}\def\elemma{\end{lemma}}}
 \def\btheorem{\begin{theorem}\sl{}\def\etheorem{\end{theorem}}}
 \def\bcorollary{\begin{corollary}\sl{}\def\ecorollary{\end{corollary}}}
 
 \def\bproposition{\begin{proposition}\sl{}\def\eproposition{\end{proposition}}}
 \def\bremark{\begin{remark}\sl{}\def\eremark{\end{remark}}}

 \def\beqlb{\begin{eqnarray}}\def\eeqlb{\end{eqnarray}}
 \def\beqnn{\begin{eqnarray*}}\def\eeqnn{\end{eqnarray*}}

 \def\<{\langle}\def\>{\rangle}

 \def\eqref#1{{\rm(\ref{#1})}}

\def\D{\textup{D}}

\def\d{\textup{d}}
\def\I{\textup{I}}
\def\e{\textup{e}}

\def\fin{\hfill$\square$}

\def\newdot{{\kern.8pt\cdot\kern.8pt}}

\def\R{\mathbb{R}}
\def\E{\mathbb{E}}
\def\P{\mathbb{P}}

\def\D{\mathbb{D}}

\def\<{\langle}
\def\>{\rangle}
\def\Proof.{\noindent{\bf Proof.}}

\begin{document}

\

\noindent{}

\bigskip\bigskip

\centerline{\Large\bf A study on the fractional Gruschin type process
\footnote{Supported by the National Natural Science Foundation of China (Grant No. 11871076),
the Natural Science Foundation of Anhui Province (Grant No. 1908085MA07).
}}

\smallskip
\
\bigskip\bigskip

\centerline{Xiliang Fan$^\dag$, Rong Yu$^\ddag$\footnote{E-mail address: fanxiliang0515@163.com (X. Fan), yurong1108@126.com(R. Yu).}}
\bigskip

\centerline{Department of Statistics, Anhui Normal University, Wuhu 241003, China}
\smallskip

\smallskip

\bigskip\bigskip

{\narrower{\narrower

\noindent{\bf Abstract.} In this article, we first establish derivative formulae for fractional Gruschin type process,
which generalize the result of Wang (J Theor Probab 27:80--95, Theorem 1.1, 2012).
Since we work on a non-Markovian context, some technical difficulties appear in the study.
Then, using the fractional calculus technique, we also derive the gradient estimate.
}


\bigskip


%%%%%%%%%%%%%%% Introduction %%%%%%%%%%%%

\section{Introduction}

\setcounter{equation}{0}

For the usual $(m+d)$-dimensional Browinian motion $(B,\tilde{B})$ and $\sigma\in C^1(\R^m;\R^d\times\R^d)$ that might be degenerate,
consider the following stochastic differential equation (SDE) on $\R^{m+d}$:
\begin{equation}\label{1.1}
\begin{cases}
 \textnormal\d X_t=\d B_t,\\
 \textnormal\d Y_t=\sigma(X_t)\d\tilde{B}_t.
\end{cases}
\end{equation}
It is easy to see that the associated generator is the following Gruschin type operator on $\R^{m+d}$:
\begin{equation*}
L(x,y)=\frac{1}{2}\left[\sum\limits_{i=1}^m\partial_{x_i}^2+\sum\limits_{j,k=1}^d(\sigma(x)\sigma(x)^*)_{jk}\partial_{y_j}\partial_{y_k}\right]
\ \ \ \ \ \mathrm{with}\ \ \ \ \ \ (x,y)\in\R^{m+d}.
\end{equation*}
If $m=d=1$ and $\sigma(x)=x$, it reduces to the Gruschin operator.
In \cite{Wang12a} Bismut type derivative formula and gradient estimates were established for the semigroup generated by the operator $L$ i.e.,
the Gruschin type semigroup.
The argument consists of using Malliavin calculus techniques together with solving a control problem.
Remark that, as pointed out by the author in \cite{Wang12a},
it seems hard to adopt the same arguments developed for Heisenberg group and subelliptic operators satisfying the generalized curvature.

In probability theory, the gradient estimate problem is an interesting research topic and has been extensively studied.
Among many probabilistic methods,
the derivative formula, which is called the Bismut formula or the Bismut-Elworthy-Li formula due to \cite{Bismut84,Elworthy&Li94},
is known to be quite effective.
As far as we are aware of, in the Markovian setting there exist many results on the derivative formulae.
Here we would like to mention a few of them.
Based upon martingale method, coupling argument, or Malliavin calculus,
the derivative formulae are derived in diffusion processes \cite{Wang11a},
degenerate SDEs \cite{Bao&Wang&Yuan13b,Guillin&Wang12,Priola06,Wang&Zhang13,Zhang10a},
and stochastic partial differential equations (SPDEs) \cite{Bao&Wang&Yuan13c,Dong&Xie10,Wang&Xu10a}.
Besides, in \cite{Wang&Xu14a} log-Harnack inequality of the Gruschin type semigroup is established by using coupling by change of measure.
Afterwards, Deng and Zhang \cite{Deng&Zhang18a} extended the results to the Markov semigroups generated by non-local Gruschin type operators
via coupling in two steps together with the regularization approximations of the underlying subordinators.

Recently, the theory of SDE driven by fractional Brownian motion was developed.
As a result, an analogous gradient estimate problem is gather attention.
The SDE above is usually understood in the sense of the pathwise Riemann-Stieltjes integration (or Young integration)
originally due to \cite{Young36a} and developed in \cite{Zahle98a}, and the rough path introduced in \cite{Lyons98a}.
See e.g. \cite{Boufoussi&Hajji11,Coutin&Qian02a,Diopa&Garrido-Atienza14,Ferrante&Rovira06,Neuenkirch&Nourdin&Tindel08,Nualart&Rascanu02a}
for the study of the existence and uniqueness of solutions to such stochastic equations.
For other interesting results involved with the distributional and paths regularities,
one may refer to \cite{Xiao11a,Baudoin&Ouyang&Tindel11a,Hairer&Pillai11,Nourdin&Simon06,Saussereau12,Yan16a} and references within.
In the previous papers \cite{Fan13b} and \cite{Fan13a,Fan19a},
we studied the derivative formulae and the Harnack type inequalities for SDEs with additive fractional noises for $H<1/2$ and $H>1/2$, respectively.

The main purpose of this article is to provide the derivative formulae and the gradient estimate for multidimensional SDE driven by fractional Brownian motion,
which is a ``fractional version" of the Gruschin type process and allows the coefficient to be degenerate.
This can be regarded as an attempt toward the solving of the derivative formulae
for more general SDEs with multiplicative fractional noise.
Although the basic strategy of the proof relies on Malliavin calculus which is similar to the case in \cite{Wang12a},
the proof gets much  more technically difficult since we work on a non-Markovian context.
We will carry it out using fractional integral and derivative operators combined with the transfer principle.

The structure of this paper is organized in the following way.
Section 2 includes some basic results about fractional calculus and fractional Brownian motion.
In Section 3, we state the derivative formulae and then present the gradient estimate
for multidimensional SDE with fractional noise which is a ``fractional version" of the Gruschin type process via Malliavin calculus.
Section 4 contains the study of the derivative formula of a more general model.

Throughout the paper, we will denote by $C$ a positive constant that may vary from one formula to another.

\section{Preliminaries}

\setcounter{equation}{0}

\subsection{Fractional calculus}

\setcounter{equation}{0}

For later use, we introduce some basic definitions and results about fractional calculus.
An exhaustive survey can be found in \cite{Samko&Kilbas&Marichev}.

Let $a,b\in\R$ with $a<b$.
For $f\in L^1(a,b)$ and $\alpha>0$,
the left-sided (resp. right-sided) fractional Riemann-Liouville integral of $f$ of order $\alpha$ on $[a,b]$ is given at almost all $x$ by
\beqnn
I_{a+}^\alpha f(x)=\frac{1}{\Gamma(\alpha)}\int_a^x\frac{f(y)}{(x-y)^{1-\alpha}}\d y,\ \
\left(\mbox{resp.}~I_{b-}^\alpha f(x)=\frac{(-1)^{-\alpha}}{\Gamma(\alpha)}\int_x^b\frac{f(y)}{(y-x)^{1-\alpha}}\d y\right),
\eeqnn
where $(-1)^{-\alpha}=\e^{-i\alpha\pi},\Gamma$ stands for the Euler function.
The above integral extends the usual $n$-order iterated integrals of $f$ for $\alpha=n\in\mathbb{N}$.

The fractional derivative may be introduced as an inverse operation.
Let $\alpha\in(0,1)$ and $p\geq1$.
If $f\in I_{a+}^\alpha(L^p([a,b],\R))$ (resp. $I_{b-}^\alpha(L^p([a,b],\R)))$, the function $\phi$ satisfying $f=I_{a+}^\alpha\phi$ (resp. $f=I_{b-}^\alpha\phi$) is unique in $L^p([a,b],\R)$
and it agrees with the left-sided (resp. right-sided) Riemann-Liouville derivative
of $f$ of order $\alpha$ given by
\beqnn
D_{a+}^\alpha f(x)=\frac{1}{\Gamma(1-\alpha)}\frac{\d}{\d x}\int_a^x\frac{f(y)}{(x-y)^\alpha}\d y\quad
\left(\mbox{resp.}\ D_{b-}^\alpha f(x)=\frac{(-1)^{1+\alpha}}{\Gamma(1-\alpha)}\frac{\d}{\d x}\int_x^b\frac{f(y)}{(y-x)^\alpha}\d y\right).
\eeqnn
The corresponding Weyl representation reads as follows
\beqnn
D_{a+}^\alpha f(x)=\frac{1}{\Gamma(1-\alpha)}\left(\frac{f(x)}{(x-a)^\alpha}+\alpha\int_a^x\frac{f(x)-f(y)}{(x-y)^{\alpha+1}}\d y\right)\\
\left(\mbox{resp.}\ \ D_{b-}^\alpha f(x)=\frac{(-1)^\alpha}{\Gamma(1-\alpha)}\left(\frac{f(x)}{(b-x)^\alpha}+\alpha\int_x^b\frac{f(x)-f(y)}{(y-x)^{\alpha+1}}\d y\right)\right),
\eeqnn
where the convergence of the integrals at the singularity $y=x$ holds pointwise for almost all $x$ if $p=1$ and in $L^p$-sense if $1<p<\infty$.

Suppose that $f\in C^\lambda(a,b)$ (the set of $\lambda$-H\"{o}lder continuous functions on $[a,b]$) and $g\in C^\mu(a,b)$ with $\lambda+\mu>1$.
The Riemann-Stieltjes integral $\int_a^bf\d g$ exists by the results of Young \cite{Young36a}.
In \cite{Zahle98a}, Z\"{a}hle presents an explicit expression for the integral $\int_a^bf\d g$ in terms of fractional derivatives.
Let $\lambda>\alpha$ and $\mu>1-\alpha$.
Then the Riemann-Stieltjes integral has the following representation:
\beqlb\label{2.0}
\int_a^bf\d g=(-1)^\alpha\int_a^bD_{a+}^\alpha f(t)D_{b-}^{1-\alpha}g_{b-}(t)\d t,\nonumber
\eeqlb
where $g_{b-}(t)=g(t)-g(b)$.
This can be regarded as fractional integration by parts formula.

\subsection{Fractional Brownian motion}

\setcounter{equation}{0}

In this part, we shall recall some important definitions and results concerning the fractional Brownian motion.
For a deeper discussion, we refer the reader to \cite{Alos&Mazet&Nualart01a,Biagini&Hu08a,Decreusefond&Ustunel98a,Mishura08a,Nualart06a}.

The $d$-dimensional fractional Brownian motion with Hurst parameter $H\in(0,1)$ on the probability space $(\Omega,\mathscr{F},\mathbb{P})$
can be defined as the centered Gauss process $B^H=\{B_t^H, t\in[0,T]\}$ with covariance function $\E\left(B_t^{H,i}B_s^{H,j}\right)=R_H(t,s)\delta_{i,j}$,
where
\beqnn
 R_H(t,s)=\frac{1}{2}\left(t^{2H}+s^{2H}-|t-s|^{2H}\right).
\eeqnn
If $H=1/2$, the process $B^H$ is a standard $d$-dimensional Brownian motion.
By the above covariance function, one can show that $\E|B_t^{H,i}-B_s^{H,i}|^p=C(p)|t-s|^{pH},\ \forall p\geq 1$.
Then, by the Kolmogorov continuity criterion
$B^{H,i}$ have $(H-\epsilon)$-order H\"{o}lder continuous paths for all $\epsilon>0,\ i=1,\cdot\cdot\cdot,d$.

For each $t\in[0,T]$, we denote by $\mathcal {F}_t$ the $\sigma$-field generated by the random variables $\{B_s^H:s\in[0,t]\}$ and the
sets of probability zero.

We denote by $\mathscr{E}$ the set of step functions on $[0,T]$.
Let $\mathcal {H}$ be the Hilbert space defined as the closure of
$\mathscr{E}$ with respect to the scalar product
\beqnn
\langle (I_{[0,t_1]},\cdot\cdot\cdot,I_{[0,t_d]}),(I_{[0,s_1]},\cdot\cdot\cdot,I_{[0,s_d]})\rangle_\mathcal {H}=\sum\limits_{i=1}^dR_H(t_i,s_i).
\eeqnn
By bounded linear transform theorem,
the mapping $(I_{[0,t_1]},\cdot\cdot\cdot,I_{[0,t_d]})\mapsto\sum_{i=1}^dB_{t_i}^{H,i}$ can be extended to an isometry between $\mathcal {H}$ and the Gaussian space $\mathcal {H}_1$ associated with $B^H$.
We denote this isometry by $\phi\mapsto B^H(\phi)$.

We recall that by \cite{Decreusefond&Ustunel98a} the covariance kernel $R_H(t,s)$ can be written as
\beqnn
 R_H(t,s)=\int_0^{t\wedge s}K_H(t,r)K_H(s,r)\d r,
\eeqnn
where $K_H$ is a square integrable kernel given by
\beqnn
K_H(t,s)=\Gamma\left(H+\frac{1}{2}\right)^{-1}(t-s)^{H-\frac{1}{2}}F\left(H-\frac{1}{2},\frac{1}{2}-H,H+\frac{1}{2},1-\frac{t}{s}\right),
\eeqnn
in which $F(\cdot,\cdot,\cdot,\cdot)$ is the Gauss hypergeometric function (for details see \cite{Decreusefond&Ustunel98a} or \cite{Nikiforov&Uvarov88}).

Now, consider the linear operator $K_H^*:\mathscr{E}\rightarrow L^2([0,T],\R^d)$ defined by
\beqnn
(K_H^*\phi)(s)=K_H(T,s)\phi(s)+\int_s^T(\phi(r)-\phi(s))\frac{\partial K_H}{\partial r}(r,s)\d r.
\eeqnn
By integration by parts, it is easy to see that when $H>1/2$, the above relation can be rewritten as
\beqlb\label{2.0'}
(K_H^*\phi)(s)=\int_s^T\phi(r)\frac{\partial K_H}{\partial r}(r,s)\d r.
\eeqlb
Due to \cite{Alos&Mazet&Nualart01a}, for any $\phi,\psi\in\mathscr{E}$,
there holds $\langle K_H^*\phi,K_H^*\psi\rangle_{L^2([0,T],\R^d)}=\langle\phi,\psi\rangle_\mathcal {H}$
and then $K_H^*$ can be extended to an isometry between $\mathcal{H}$ and $L^2([0,T],\R^d)$.
So, owing to \cite{Alos&Mazet&Nualart01a} again,
the process $\{W_t=B^H((K_H^*)^{-1}{\rm I}_{[0,t]}),t\in[0,T]\}$ is a Wiener process,
and $B^H$ has Volterra's representation of the form
\beqlb\label{2.0''}
 B^H_t=\int_0^tK_H(t,s)\d W_s.
\eeqlb
On the other hand, define the operator $K_H: L^2([0,T],\mathbb{R}^d)\rightarrow I_{0+}^{H+1/2 }(L^2([0,T],\mathbb{R}^d))$
associated with the kernel $K_H(\cdot,\cdot)$ by
\beqnn
 (K_Hf^i)(t)=\int_0^tK_H(t,s)f^i(s)\d s,\ \ i=1,\cdot\cdot\cdot,d.
\eeqnn
By \cite{Decreusefond&Ustunel98a}, it is an isomorphism and can be expressed in terms of fractional integrals as follows:
\beqnn
 (K_H f)(s)=I_{0+}^{1}s^{H-1/2}I_{0+}^{H-1/2}s^{1/2-H}f,\ H\geq1/2,\\
 (K_H f)(s)=I_{0+}^{2H}s^{1/2-H}I_{0+}^{1/2-H}s^{H-1/2}f,\ H\leq1/2.
\eeqnn
where $f\in L^2([0,T],\mathbb{R}^d)$.
Consequently, for every $h\in I_{0+}^{H+1/2}(L^2([0,T],\R^d))$, the inverse operator $K_H^{-1}$
is of the following form
\beqlb\label{2.1}
(K_H^{-1}h)(s)=s^{H-1/2}D_{0+}^{H-1/2}s^{1/2-H}h',\ H\geq1/2,
\eeqlb
\beqlb\label{2.2}
(K_H^{-1}h)(s)=s^{1/2-H}D_{0+}^{1/2-H}s^{H-1/2}D_{0+}^{2H}h,\ H\leq1/2.
\eeqlb

The remaining part will be devoted to the Malliavin calculus of fractional Brownian motion.

Let $\Omega$ be the canonical probability space $C_0([0,T],\R^d)$, the set of continuous functions,
null at time $0$, equipped with the supremum norm.
Let $\P$ be the unique probability measure on $\Omega$
such that the canonical process $\{B^H_t; t\in[0,T]\}$ is a $d$-dimensional fractional Brownian motion with Hurst parameter $H$.
Then, the injection $R_H=K_H\circ K_H^*:\mathcal{H}\rightarrow\Omega$ embeds $\mathcal{H}$ densely into $\Omega$ and
$(\Omega,\mathcal{H},\P)$ is an abstract Wiener space in the sense of Gross.
In the sequel we will make this assumption on the underlying probability space.

We denote by $\mathcal {S}$ the set of smooth and cylindrical random variables of the form
\beqnn
F=f\left(B^H(\phi_1),\cdot\cdot\cdot,B^H(\phi_n)\right),
\eeqnn
where $n\geq 1, f\in C_b^\infty(\mathbb{R}^n)$, the set of $f$ and all its partial derivatives are bounded, $\phi_i\in\mathcal{H}, 1\leq i\leq n$.
The Malliavin derivative of $F$, denoted by $\mathbb{D}F$, is defined as the $\mathcal {H}$-valued random variable
\beqnn
\mathbb{D}F=\sum_{i=1}^n\frac{\partial f}{\partial x_i}\left(B^H(\phi_1),\cdot\cdot\cdot,B^H(\phi_n)\right)\phi_i.
\eeqnn
For any $p\geq 1$, we define the Sobolev space $\mathbb{D}^{1,p}$ as the completion of $\mathcal {S}$ w.r.t. the norm
\beqnn
\|F\|_{1,p}^p=\mathbb{E}|F|^p+\mathbb{E}\|\mathbb{D}F\|^p_{\mathcal {H}}.
\eeqnn
While we will denote by $\delta$ and $\mathrm{Dom}\delta$ the divergence operator of $\mathbb{D}$ and its domain.
We conclude this section by giving a transfer principle that connects the derivative and divergence operators of both processes $B^H$ and $W$.

\bproposition\label{Pro 2.1}\cite[Proposition 5.2.1]{Nualart06a}
For any $F\in\mathbb{D}_W^{1,2}=\mathbb{D}^{1,2}$, there holds
\beqnn
K_H^*\mathbb{D}F=\mathbb{D}_WF,
\eeqnn
where $\mathbb{D}_W$ denotes the derivative operator w.r.t. $W$,
and $\mathbb{D}_W^{1,2}$ the corresponding Sobolev space.
\eproposition

\bproposition\label{Pro 2.2}\cite[Proposition 5.2.2]{Nualart06a}
$\mathrm{Dom}\delta=(K_H^*)^{-1}(\mathrm{Dom}\delta_W)$,
and for any $\mathcal{H}$-valued random variable $u$ in $\mathrm{Dom}\delta$ we have $\delta(u)=\delta_W(K_H^*u)$,
where $\delta_W$ denotes the divergence operator w.r.t. $W$.
\eproposition

\section{Derivative formulae and gradient estimate}

\setcounter{equation}{0}

The objective of this section is to study the following SDE on $\R^{d_1+d_2}$:
\begin{equation}\label{3.1}
\begin{cases}
 \textnormal\d X_t=\d B_t^H,\\
 \textnormal\d Y_t=\sigma(X_t)\d\tilde{B}_t^H,
\end{cases}
\end{equation}
where $(B_t^H,\tilde{B}_t^H)_{t\in[0,T]}$ is a fractional Brownian motion on  $\R^{d_1+l}$ with $H>1/2$,
$\sigma:\R^{d_1}\rightarrow\R^{d_2}\times\R^{l}$.

Remark that, the coefficient $\sigma$ is allowed to be degenerate,
and compared with \eqref{1.1}, the equation \eqref{3.1} can be viewed as a ``fractional version" of the Gruschin type process.
We will use $(X_t,Y_t)$ to denote the solution and denote $\E^{x,y}$ by the expectation taken for the solution  with the initial value $(x,y)\in\R^{d_1+d_2}$.
In this part, we aim to establish the Bismut type derivative formulae for the associated family of operators $(P_t)_{0\leq t\leq T}$:
$$P_tf(x,y)=\E^{x,y}f(X_t,Y_t),\ \ \ (x,y)\in\R^{d_1+d_2},f\in\mathcal{B}_b(\R^{d_1+d_2}),$$
where $\mathcal{B}_b(\R^{d_1+d_2})$ is the set of all bounded measurable functions on $\R^{d_1+d_2}.$
Besides, for $f\in C^\alpha([0,T];\R^m)$, set
\beqnn
\|f\|_\infty:=\sup\limits_{t\in[0,T]}|f(t)|,\ \ \ \|f\|_\alpha:=\sup\limits_{s\neq t,s,t\in[0,T]}\frac{|f(t)-f(s)|}{|t-s|^\alpha}.
\eeqnn

To the end, we first consider the following more general SDE with fractional noise on the time interval $[0,T]$:
\beqlb\label{3.2}
\d Z_t=\bar{b}(Z_t)\d t+\bar{\sigma}(Z_t)\d \bar{B}^H(t), \ Z_0=z\in\R^d,
\eeqlb
where $\bar{b}:\R^d\rightarrow\R^d,\bar{\sigma}:\R^d\rightarrow\R^{d}\times\R^m,\bar{B}^H$ is a $m$-dimensional fractional Brownian motion with $H>1/2$,
and then introduce the assumption (H1): $\bar{b}$ and $\bar{\sigma}$ are both differentiable with bounded derivatives.

Below we shall present a lemma which is important for the proof our main results.
\blemma\label{Lem 3.1}
Assume (H1).
Then, there is a unique solution $Z$ of \eqref{3.2} with $Z^i(t)\in\D^{1,2},i=1,\cdots,d, t\in[0,T]$.
Moreover, there holds, for any $h\in\mathcal {H}$,
\beqnn
\mathbb{D}_{R_Hh}Z_t=\int_0^t\bar{\sigma}(Z_u)\d (R_H h)(u)
+\int_0^t\nabla \bar{b}(Z_u)\mathbb{D}_{R_Hh}Z_u\d u
+\int_0^t\nabla\bar{\sigma}(Z_u)\mathbb{D}_{R_Hh}Z_u\d \bar{B}_u.
\eeqnn
\elemma
\emph{Proof.}
Owing to \cite[Theorem 2.1]{Nualart&Rascanu02a},
the condition (H1) ensures that there exists a unique adapted solution $Z$ to the equation \eqref{3.2}.
Moreover, by \cite[Theorem 6 and Proposition 7]{Nualart&Saussereau09}, it follows that $Z_t^i$ belongs to $\D^{1,2}$ and satisfies,
for $i=1,\cdots,d, j=1,\cdots,m$,
\begin{equation}\label{Pf of Lem 3.1-1}
\mathbb{D}_s^jZ^i_t=
\left\{
\begin{array}{ll}
\bar{\sigma}_{ij}(Z_s)+\sum\limits_{k=1}^d\int_s^t\partial_k \bar{b}_i(Z_u)\mathbb{D}_s^jZ^k_u\d u
+\sum\limits_{k=1}^d\sum\limits_{l=1}^m\int_s^t\partial_k\bar{\sigma}_{il}(Z_u)\mathbb{D}_s^jZ^k_u\d \bar{B}_u^l, \ s\leq t,\\
0, \  {\rm else}.
\end{array} \right.
\end{equation}
We can reformulate \eqref{Pf of Lem 3.1-1} as
\beqlb\label{Pf of Lem 3.1-2}
\mathbb{D}Z^i_t&=&\sum\limits_{l=1}^m(\bar{\sigma}_{il}(Z)\I_{[0,t]})e^l
+\sum\limits_{k=1}^d\int_0^t\partial_k \bar{b}_i(Z_u)\mathbb{D}Z^k_u\d u\nonumber\\
&&+\sum\limits_{k=1}^d\sum\limits_{l=1}^m\int_0^t\partial_k\bar{\sigma}_{il}(Z_u)\mathbb{D}Z^k_u\d \bar{B}_u^l,
\ \ \ \ \ 1\leq i\leq d, \ 0\leq t\leq T,\nonumber
\eeqlb
where $\{e^l\}_{l=1}^m$ is the canonical ONB on $\R^m$.\\
Then, we have, for each $h\in\mathcal {H}$,
\beqlb\label{Pf of Lem 3.1-2'}
\langle\mathbb{D}Z^i_t,h\rangle_\mathcal{H}&=&\sum\limits_{l=1}^m\langle(\bar{\sigma}_{il}(Z)\I_{[0,t]})e^l,h\rangle_\mathcal{H}
+\sum\limits_{k=1}^d\int_0^t\partial_k \bar{b}_i(Z_u)\langle\mathbb{D}Z^k_u,h\rangle_\mathcal{H}\d u\nonumber\\
&&+\sum\limits_{k=1}^d\sum\limits_{l=1}^m\int_0^t\partial_k\bar{\sigma}_{il}(Z_u)\langle\mathbb{D}Z^k_u,h\rangle_\mathcal{H}\d \bar{B}_u^l.
\eeqlb
Observe that, for each $i=1,\cdots,d,l=1,\cdots,m$, by \eqref{2.0'} we get
\beqlb\label{Pf of Lem 3.1-3}
\left\langle(\bar{\sigma}_{il}(Z)\I_{[0,t]})e^l,h\right\rangle_\mathcal {H}
&=&\langle K_H^*((\bar{\sigma}_{il}(Z)\I_{[0,t]})e^l),K_H^*h\rangle_{L^2([0,T],\R^d)}\nonumber\\
&=&\sum_{j=1}^d\int_0^T\left(K_H^*((\bar{\sigma}_{il}(Z)\I_{[0,t]})e^l)\right)^j(s)(K_H^*h)^j(s)\d s\nonumber\\
&=&\int_0^T\left(K_H^*(\bar{\sigma}_{il}(Z)\I_{[0,t]})\right)(s)(K_H^*h)^l(s)\d s\nonumber\\
&=&\int_0^T\int_s^T\bar{\sigma}_{il}(Z_r)\I_{[0,t]}(r)\frac{\partial K_H}{\partial r}(r,s)(K_H^*h)^l(s)\d r\d s\nonumber\\
&=&\int_0^t\bar{\sigma}_{il}(Z_r)\int_0^r\frac{\partial K_H}{\partial r}(r,s)(K_H^*h)^l(s)\d s\d r\nonumber\\
&=&\int_0^t\bar{\sigma}_{il}(Z_r)\d (R_H h)^l(r),
\eeqlb
where the last equality is due to the fact: $R_H=K_H\circ K_H^\ast$.

Next, by \cite[page 400]{Nualart&Saussereau09} or
following the arguments of \cite[Lemma 3.1 and Proposition 3.3]{Fan13a}, we deduce that
$\mathbb{D}_{R_Hh}Z^i_t=\frac{\d}{\d\epsilon}\Big|_{\epsilon=0}Z^i_t(w+\epsilon R_Hh),\ h\in\mathcal {H},$ satisfies
\beqlb\label{Pf of Lem 3.1-4}
\mathbb{D}_{R_Hh}Z^i_t=\langle\mathbb{D}Z^i_t,h\rangle_\mathcal{H},\ \ \ 1\leq i\leq d.
\eeqlb
So, plugging \eqref{Pf of Lem 3.1-3} and \eqref{Pf of Lem 3.1-4} into \eqref{Pf of Lem 3.1-2'} yields
\beqlb\label{Pf of Lem 3.1-5}
\mathbb{D}_{R_Hh}Z^i_t&=&\sum\limits_{l=1}^m\int_0^t\bar{\sigma}_{il}(Z_u)\d (R_H h)^l(u)
+\sum\limits_{k=1}^d\int_0^t\partial_k \bar{b}_i(Z_u)\mathbb{D}_{R_Hh}Z^k_u\d u\nonumber\\
&&+\sum\limits_{k=1}^d\sum\limits_{l=1}^m\int_0^t\partial_k\bar{\sigma}_{il}(Z_u)\mathbb{D}_{R_Hh}Z^k_u\d \bar{B}_u^l.\nonumber
\eeqlb
Then we obtain
\beqlb\label{Pf of Lem 3.1-6}
\mathbb{D}_{R_Hh}Z_t&:=&
\left(
\begin{array}{ccc}
\mathbb{D}_{R_Hh}Z^1_t\\
\cdot\\
\cdot\\
\cdot\\
\mathbb{D}_{R_Hh}Z^d_t
\end{array}
\right)\nonumber\\
&=&\int_0^t\bar{\sigma}(Z_u)\d (R_H h)(u)
+\int_0^t\nabla \bar{b}(Z_u)\mathbb{D}_{R_Hh}Z_u\d u
+\int_0^t\nabla\bar{\sigma}(Z_u)\mathbb{D}_{R_Hh}Z_u\d \bar{B}_u.\nonumber
\eeqlb
The proof is completed.
\fin

Now, we turn to the equation \eqref{3.1}.
For any $v=(v_1,v_2)\in\R^{d_1}\times\R^{d_2}$, we are to find $h=(h_1,h_2)\in {\rm Dom}\delta$ such that
\begin{equation}\label{3.3}
\nabla_v P_Tf(x,y)=\E^{x,y}(f(X_T,Y_T)\delta(h)),\ f\in C_b^1(\R^{d_1+d_2})
\end{equation}
holds.
We will show that $h$ satisfies \eqref{3.3} on condition that  it is in ${\rm Dom}\delta$ with
\begin{equation}\label{3.4'}
(R_Hh_1)(0)=0,\ (R_Hh_1)(T)=v_1
\end{equation}
and
\begin{equation}\label{3.4}
\int_0^T\sigma(X_u)(R_H h_2)'(u)\d u+\int_0^T\nabla\sigma(X_u)((R_Hh_1)(u)-v_1)\d\tilde{B}^H_u-v_2=0.
\end{equation}

\bproposition\label{Pro 3.1}
Suppose that $\sigma$ is differentiable with bounded derivative.
For $v=(v_1,v_2)\in\R^{d_1+d_2}$, let $h=(h_1,h_2)$ be given in \eqref{3.4'} and  \eqref{3.4}.
If $h\in{\rm Dom}\delta$, then there holds \eqref{3.3}.
\eproposition
\emph{Proof.}
According to Lemma \ref{Lem 3.1} for $d=m=d_1+d_2,Z_t=(X_t,Y_t)^*,\bar{B}^H=(B^H,\tilde{B}^H)^*,\bar{b}=0$ and
\begin{equation}\nonumber
\bar{\sigma}(x,y)=
\left(
\begin{array}{ccc}
I_{d_1\times d_1} & 0\\
0 & \sigma(x)\\
\end{array}
\right)
\end{equation}
being in \eqref{3.2}, we have, for each $h=(h_1,h_2)\in\mathcal {H}$,
\begin{equation}\label{Pf of Lem 3.2-1}\nonumber
\left\{
\begin{array}{ll}
\mathbb{D}_{R_Hh}X_t=(R_Hh_1)(t)-(R_Hh_1)(0),\\
\mathbb{D}_{R_Hh}Y_t=\int_0^t\sigma(X_u)(R_H h_2)'(u)\d u
+\int_0^t\nabla\sigma(X_u)\mathbb{D}_{R_Hh}X_u\d\tilde{B}^H_u.
\end{array} \right.
\end{equation}
This, together with $(R_Hh_1)(0)=0$, leads to
\begin{equation}\label{Pf of Lem 3.2-2}
\left\{
\begin{array}{ll}
\mathbb{D}_{R_Hh}X_t=(R_Hh_1)(t),\\
\mathbb{D}_{R_Hh}Y_t=\int_0^t\sigma(X_u)(R_H h_2)'(u)\d u
+\int_0^t\nabla\sigma(X_u)(R_Hh_1)(u)\d\tilde{B}^H_u.
\end{array} \right.
\end{equation}
On the other hand, we easily know that the directional derivative processes satisfy that,
for any $(v_1,v_2)\in\R^{d_1}\times\R^{d_2}$,
\begin{equation}\label{Pf of Lem 3.2-3}\nonumber
\left\{
\begin{array}{ll}
\nabla_vX_t=v_1,\\
\nabla_vY_t=v_2+\int_0^t\nabla\sigma(X_u)\nabla_vX_u\d\tilde{B}^H_u.
\end{array} \right.
\end{equation}
Hence, we get
\begin{equation}\label{Pf of Lem 3.2-4}
\left\{
\begin{array}{ll}
\nabla_vX_t=v_1,\\
\nabla_vY_t=v_2+\int_0^t\nabla\sigma(X_u)v_1\d\tilde{B}^H_u.
\end{array} \right.
\end{equation}
In view of the control condition \eqref{3.4} and $(R_Hh_1)(T)=v_1$, \eqref{Pf of Lem 3.2-2} and \eqref{Pf of Lem 3.2-4} imply
\begin{equation}\label{Pf of Lem 3.2-5}
(\mathbb{D}_{R_Hh}X_T,\mathbb{D}_{R_Hh}Y_T)=(\nabla_vX_T,\nabla_vY_T).\nonumber
\end{equation}
Consequently, for any $f\in C_b^1(\R^{d_1+d_2})$, we conclude
\beqlb\label{Pf of Lem 3.2-6}
\nabla_v P_T f(x,y)&=&\E^{x,y} \nabla_vf(X_T,Y_T)=\E^{x,y}((\nabla f)(X_T,Y_T)(\nabla_v X_T,\nabla_v Y_T))\cr
&=&\E^{x,y}((\nabla f)(X_T,Y_T)(\mathbb{D}_{R_Hh}X_T,\mathbb{D}_{R_Hh}Y_T))=\E^{x,y}\mathbb{D}_{R_Hh}f(X_T,Y_T)\cr
&=&\E^{x,y}\langle\mathbb{D}f(X_T,Y_T),h\rangle_\mathcal {H}=\E^{x,y}(f(X_T,Y_T)\delta(h)).
\eeqlb
\fin

With the above proposition in hand,
to establish explicit derivative formula,
we need to calculate $\delta(h)$ for $h$ which will be given a concrete choice satisfying \eqref{3.4'} and \eqref{3.4}.
To this end, we assume (H2):\\
$\sigma$ is differentiable with bounded derivative, and for any $x\in\R^{d_1}, \int_0^T(\sigma\sigma^*)\left(x+B^H_u\right)\d u$
is invertible such that
\beqnn
\E\left\|\left(\int_0^T(\sigma\sigma^*)\left(x+B^H_u\right)\d u\right)^{-1}\right\|^{2+\epsilon_0}<\infty,
\eeqnn
where $\epsilon_0>0$ is a constant.
\btheorem\label{The 3.1}
Assume (H2) and let $v=(v_1,v_2)\in\R^{d_1+d_2}$.
Then
\beqlb\label{Th3.1-0}
\nabla_v P_T f(x,y)=\E^{x,y}\left[f(X_T,Y_T)M_T\right],\ \ f\in C_b^1(\R^{d_1+d_2}),
\eeqlb
holds for
\beqlb\label{Th3.1-0'}
M_T
&=&\int_0^T\left\langle K_H^{-1}\left(\frac{\cdot}{T}\right)(t)v_1,\d W_t\right\rangle
+\left\langle\vartheta(T),\int_0^T\left(K_H^{-1}\left(\int_0^\cdot\sigma^*\left(x+B^H_u\right)\d u\right)\right)^*(t)\d\tilde{W}_t\right\rangle\cr
&&-Tr\left(\left(\int_0^T(\sigma\sigma^*)\left(x+B^H_u\right)\d u\right)^{-1}\int_0^T\frac{T-u}{T}\nabla\sigma\left(x+B^H_u\right)v_1\sigma^*(x+B^H_u)\d u\right)\nonumber\\
\eeqlb
and
\beqlb\label{Th3.1-0''}
\vartheta(T)=\left(\int_0^T(\sigma\sigma^*)\left(x+B^H_u\right)\d u\right)^{-1}\left(v_2+\int_0^T\frac{T-u}{T}\nabla\sigma\left(x+B^H_u\right)v_1\d \tilde{B}^H_u\right),
\eeqlb
where $\bar{W}:=(W,\tilde{W})$ is the underlying Wiener process w.r.t. $(B,\tilde{B})$ defined in \eqref{2.0''}.
\etheorem
\emph{Proof.}
We first take $h=(h_1,h_2)$ as follows: for $t\in[0,T]$,
\begin{equation}\label{Pf of The 3.1-1}
h_1(t)=R_H^{-1}\left(\frac{\cdot}{T}v_1\right)(t)
\end{equation}
and
\begin{equation}\label{Pf of The 3.1-2}
h_2(t)
=R_H^{-1}\left(\left(\int_0^\cdot\sigma^*(X_u)\d u\right)\left(\int_0^T(\sigma\sigma^*)(X_u)\d u\right)^{-1}\left(v_2+\int_0^T\frac{T-u}{T}\nabla\sigma(X_u)v_1\d \tilde{B}^H_u\right)\right)(t).
\end{equation}
It is easy to check that the above $h=(h_1,h_2)$ satisfies \eqref{3.4'} and \eqref{3.4}.
Now, to prove \eqref{Th3.1-0}-\eqref{Th3.1-0''},
according to Proposition \ref{Pro 3.1} it remains to verify $h=(h_1,h_2)\in{\rm Dom}\delta$ and then calculate $\delta(h)$
which will be fulfilled via Proposition \ref{Pro 2.1} and Proposition \ref{Pro 2.2} as well as \cite[Proposition 1.3.3]{Nualart06a}.
Notice that, by \eqref{Pf of The 3.1-2} we know that
\beqlb\label{Pf of The 3.1-3}
(R_Hh_2)(t)
&=&\left(\int_0^t\sigma^*(X_u)\d u\right)\left(\int_0^T(\sigma\sigma^*)(X_u)\d u\right)^{-1}\left(v_2+\int_0^T\frac{T-u}{T}\nabla\sigma(X_u)v_1\d \tilde{B}^H_u\right)\nonumber\\
&=:&\int_0^t\sigma^*(X_u)\d u\cdot\vartheta(T)
\eeqlb
is not adapted and then $(K_H^*h_2)(t)=K_H^{-1}\left(\int_0^\cdot\sigma^*(X_u)\d u\cdot\vartheta(T)\right)(t)$ is so.
To circumvent the impediment, let $\{e^i\}_{i=1}^{d_2}$ be the canonical ONB on $\R^{d_2}$, then we have
\beqlb\label{Pf of The 3.1-4}
(R_Hh_2)(t)=\int_0^t\sigma^*(X_u)\d u\left(\sum\limits_{i=1}^{d_2}\langle\vartheta(T),e^i\rangle e^i\right)
=\sum\limits_{i=1}^{d_2}\langle\vartheta(T),e^i\rangle\int_0^t\sigma^*(X_u)e^i\d u.\nonumber
\eeqlb
So, combining this with \eqref{Pf of The 3.1-1} yields
\beqlb\label{Pf of The 3.1-5}
(R_Hh)(t)&=&((R_Hh_1)(t),(R_Hh_2)(t))
=((R_Hh_1)(t),0)+(0,(R_Hh_2)(t))\nonumber\\
&=&\left(\frac{t}{T}v_1,0\right)+\sum\limits_{i=1}^{d_2}\langle\vartheta(T),e^i\rangle\left(0,\int_0^t\sigma^*(X_u)e^i\d u\right).\nonumber
\eeqlb
As a consequence, we obtain
\beqlb\label{Pf of The 3.1-6}
(K_H^*h)(t)&=&((K_H^*{h_1})(t),0)+(0,(K_H^*{h_2})(t))\nonumber\\
&=&\left(K_H^{-1}\left(\frac{\cdot}{T}\right)(t)v_1,0\right)
+\sum\limits_{i=1}^{d_2}\langle\vartheta(T),e^i\rangle\left(0,K_H^{-1}\left(\int_0^\cdot\sigma^*(X_u)e^i\d u\right)(t)\right).
\eeqlb
Obviously, we get
\beqlb\label{Pf of The 3.1-7}
\delta_{\bar{W}}(K_H^*{h_1},0)&=&\delta_{\bar{W}}\left(K_H^{-1}\left(\frac{\cdot}{T}\right)v_1,0\right)\nonumber\\
&=&\int_0^T\left\langle K_H^{-1}\left(\frac{\cdot}{T}\right)(t)v_1,\d W_t\right\rangle.
\eeqlb

Next, we shall focus on the second term of the right hand of \eqref{Pf of The 3.1-6}.
By \eqref{2.1}, we first obtain
\beqlb\label{Pf of The 3.1-8}
&&K_H^{-1}\left(\int_0^\cdot\sigma^*(X_u)e^i\d u\right)(t)
=t^{H-\frac{1}{2}}D_{0+}^{H-\frac{1}{2}}\left(\cdot^{\frac{1}{2}-H}\sigma^*(X_\cdot)e^i\right)(t)\cr
&=&\frac{1}{\Gamma(\frac{3}{2}-H)}\Bigg[t^{\frac{1}{2}-H}\sigma^*(X_t)e^i
+\left(H-\frac{1}{2}\right)\sigma^*(X_t)e^i t^{H-\frac{1}{2}}\int_0^t\frac{t^{\frac{1}{2}-H}-r^{\frac{1}{2}-H}}{(t-r)^{\frac{1}{2}+H}}\d r\cr
&&+\left(H-\frac{1}{2}\right)t^{H-\frac{1}{2}}
\int_0^t\frac{\sigma^*(X_t)e^i-\sigma^*(X_r)e^i}{(t-r)^{\frac{1}{2}+H}}r^{\frac{1}{2}-H}\d r\Bigg]\cr
&=:&\frac{1}{\Gamma(\frac{3}{2}-H)}[\alpha_1(t)+\alpha_2(t)+\alpha_3(t)].
\eeqlb
Using the estimate
\beqnn
|\sigma^*(X_t)e^i|=|\sigma^*(x+B^H_t)e^i|\leq\|\sigma(x)\|+C|B^H_t|
\eeqnn
and the relation
\begin{equation}\label{Pf of The 3.1-9'}
\int_0^t\frac{r^{\frac{1}{2}-H}-t^{\frac{1}{2}-H}}{(t-r)^{\frac{1}{2}+H}}\d r
=\int_0^1\frac{\theta^{\frac{1}{2}-H}-1}{(1-\theta)^{\frac{1}{2}+H}}\d\theta\cdot t^{1-2H}=:C_0t^{1-2H}
\end{equation}
with some positive constant $C_0$, we get
\beqlb\label{Pf of The 3.1-9}
|\alpha_1(t)+\alpha_2(t)|\leq Ct^{\frac{1}{2}-H}(\|\sigma(x)\|+|B^H_t|)\leq Ct^{\frac{1}{2}-H}\left(1+\|B^H\|_\infty\right).
\eeqlb
As for $\alpha_3(t)$, from (H2) it follows that
\beqlb\label{Pf of The 3.1-10}
|\alpha_3(t)|&\leq& Ct^{H-\frac{1}{2}}
\int_0^t\frac{|B^H_t-B^H_r|}{(t-r)^{\frac{1}{2}+H}}r^{\frac{1}{2}-H}\d r\cr
&\leq&Ct^{\frac{1}{2}-\epsilon}\|B^H\|_{H-\epsilon},
\eeqlb
where $\epsilon$ is fixed satisfying $0<\epsilon<\frac{1}{2}$.\\
Then, by \eqref{Pf of The 3.1-8}-\eqref{Pf of The 3.1-10} and the Fernique theorem (see, for instance, \cite{Berman85a}),
we derive that $K_H^{-1}\left(\int_0^\cdot\sigma^*(X_u)e^i\d u\right)\in L^2([0,T]\times\Omega,\R^{d_2})$.
Moreover, it is clear that
$K_H^{-1}\left(\int_0^\cdot\sigma^*(X_u)e^i\d u\right)$ is adapted due to the fact that the operator $K_H^{-1}$ preserves the adaptability property.
So, we have
\beqlb\label{Pf of The 3.1-11}
\delta_{\bar{W}}\left(0,K_H^{-1}\left(\int_0^\cdot\sigma^*(X_u)e^i\d u\right)\right)
=\int_0^T\left\langle K_H^{-1}\left(\int_0^\cdot\sigma^*(X_u)e^i\d u\right)(t),\d \tilde{W}_t\right\rangle.
\eeqlb
Now, let $\mathcal{F}_T=\sigma(B_t^H:0\leq t\leq T)$.
By \eqref{Pf of The 3.1-11}, the definition of $\vartheta(T)$, the $C_r$-inequality and the H\"{o}lder inequality, we obtain
\beqlb\label{Pf of The 3.1-12}
&&\E\left(\langle\vartheta(T),e^i\rangle\cdot\delta_{\bar{W}}\left(0,K_H^{-1}\left(\int_0^\cdot\sigma^*(X_u)e^i\d u\right)\right)\right)^2\cr
&=&
\E\left\{\E\left[\left(\langle\vartheta(T),e^i\rangle\cdot\int_0^T
\left\langle K_H^{-1}\left(\int_0^\cdot\sigma^*(X_u)e^i\d u\right)(t),\d \tilde{W}_t\right\rangle\right)^2\Bigg|\mathcal{F}_T\right]\right\}\cr
&\leq&
2\E\left\{\left|\left\langle\left(\int_0^T(\sigma\sigma^*)(X_u)\d u\right)^{-1}v_2,e^i\right\rangle\right|^2\cdot\E\left[\left(\int_0^T
\left\langle K_H^{-1}\left(\int_0^\cdot\sigma^*(X_u)e^i\d u\right)(t),\d \tilde{W}_t\right\rangle\right)^2\Bigg|\mathcal{F}_T\right]\right\}\cr
&&
+2\E\Bigg\{\left[\E\left|\left\langle\left(\int_0^T(\sigma\sigma^*)(X_u)\d u\right)^{-1}\int_0^T\frac{T-u}{T}\nabla\sigma(X_u)v_1\d \tilde{B}^H_u,e^i\right\rangle\right|^4\Bigg|\mathcal{F}_T\right]^{\frac{1}{2}}\cr
&&\ \ \ \ \ \  \times\left[\E
\left|\int_0^T\left\langle K_H^{-1}\left(\int_0^\cdot\sigma^*(X_u)e^i\d u\right)(t),\d \tilde{W}_t\right\rangle\right|^4\Bigg|\mathcal{F}_T\right]^{\frac{1}{2}}\Bigg\}\cr
&\leq&
C\E\left\{\left\|\left(\int_0^T(\sigma\sigma^*)(X_u)\d u\right)^{-1}\right\|^2\cdot\int_0^T\left|K_H^{-1}\left(\int_0^\cdot\sigma^*(X_u)e^i\d u\right)(t)\right|^2\d t\right\}\cr
&&
+C\E\Bigg\{\left[\E\left|\int_0^T\frac{T-u}{T}\left\langle\left(\left(\int_0^T(\sigma\sigma^*)(X_u)\d u\right)^{-1}\nabla\sigma(X_u)v_1\right)^*e^i,\d\tilde{B}^H_u\right\rangle\right|^4\Bigg|\mathcal{F}_T\right]^{\frac{1}{2}}\cr
&&\ \ \ \ \ \ \ \times\left[\int_0^T\left|K_H^{-1}\left(\int_0^\cdot\sigma^*(X_u)e^i\d u\right)(t)\right|^4\d t\right]^{\frac{1}{2}}\Bigg\}\cr
&\leq&
C\E\left\{\left\|\left(\int_0^T(\sigma\sigma^*)(X_u)\d u\right)^{-1}\right\|^2\cdot\int_0^T\left|K_H^{-1}\left(\int_0^\cdot\sigma^*(X_u)e^i\d u\right)(t)\right|^2\d t\right\}\cr
&&
+C\E\Bigg\{\left[\int_0^T\left|\left(\left(\int_0^T(\sigma\sigma^*)(X_u)\d u\right)^{-1}\nabla\sigma(X_u)v_1\right)^*e^i\right|^4\d u\right]^{\frac{1}{2}}\cr
&&\ \ \ \ \ \ \ \times\left[\int_0^T\left|K_H^{-1}\left(\int_0^\cdot\sigma^*(X_u)e^i\d u\right)(t)\right|^4\d t\right]^{\frac{1}{2}}\Bigg\},
\eeqlb
where the last inequality is due to \cite[page 292-293]{Nualart06a}.\\
Then, by \eqref{Pf of The 3.1-8}-\eqref{Pf of The 3.1-10} again and (H2), we deduce
\beqlb\label{Pf of The 3.1-13}
&&\E\left(\langle\vartheta(T),e^i\rangle\cdot\delta_{\bar{W}}\left(0,K_H^{-1}\left(\int_0^\cdot\sigma^*(X_u)e^i\d u\right)\right)\right)^2\cr
&\leq&
C\E\Bigg\{\left\|\left(\int_0^T(\sigma\sigma^*)(X_u)\d u\right)^{-1}\right\|^2
\Bigg(\int_0^T\left|K_H^{-1}\left(\int_0^\cdot\sigma^*(X_u)e^i\d u\right)(t)\right|^2\d t\cr
&&\ \ \ \ \ \ \ +\int_0^T\left|K_H^{-1}\left(\int_0^\cdot\sigma^*(X_u)e^i\d u\right)(t)\right|^4\d t
+\int_0^T\left\|\nabla\sigma(X_u)\right\|^4\d t\Bigg)\Bigg\}\cr
&\leq&
C\E\left\{\left\|\left(\int_0^T(\sigma\sigma^*)(X_u)\d u\right)^{-1}\right\|^2
\left(1+\|B^H\|_\infty^4+\|B^H\|_{H-\epsilon}^4\right)\right\}\cr
&\leq&C\left[\E\left\|\left(\int_0^T(\sigma\sigma^*)(X_u)\d u\right)^{-1}\right\|^{2+\epsilon_0}\right]^\frac{2}{2+\epsilon_0}\cdot
\left[\E\left(1+\|B^H\|_\infty^4+\|B^H\|_{H-\epsilon}^4\right)^{\frac{2+\epsilon_0}{\epsilon_0}}\right]^{\frac{\epsilon_0}{2+\epsilon_0}}\cr
&<&\infty.
\eeqlb
That is, $\langle\vartheta(T),e^i\rangle\cdot\delta_{\bar{W}}\left(0,K_H^{-1}\left(\int_0^\cdot\sigma^*(X_u)e^i\d u\right)\right)\in L^2(\P), \ i=1,\cdots,d_2$.\\
On the other hand, by \cite[Example 1.2.1]{Nualart06a} we know that in the case of a classical Brownian motion,
$K_H^*$ is the identity map on $L^2([0,T],\R^{d_1+d_2})$ which means $\mathcal {H}_{\bar{W}}:=\mathcal {H}=L^2([0,T],\R^{d_1+d_2})$,
and $K_H(t,s)$ equal to $I_{[0,t]}(s)$ which implies that $\mathcal {H}_H$ is the space of absolutely continuous functions,
vanishing at zero, with a square integrable derivative.
Then, using Proposition \ref{Pro 2.1} and the fact that $K_H^*$ is an isometry between $\mathcal{H}$ and $L^2([0,T],\R^{d_1+d_2})$, we obtain
\beqlb\label{Pf of The 3.1-14}
&&\left\langle\mathbb{D}_{\bar{W}}\langle\vartheta(T),e^i\rangle,\left(0,K_H^{-1}\left(\int_0^\cdot\sigma^*(X_u)e^i\d u\right)\right)\right\rangle_{\mathcal {H}_{\bar{W}}}\cr
&=&\left\langle\mathbb{D}_{\bar{W}}\langle\vartheta(T),e^i\rangle,\left(0,K_H^{-1}\left(\int_0^\cdot\sigma^*(X_u)e^i\d u\right)\right)\right\rangle_{L^2([0,T],\R^{d_1+d_2})}\cr
&=&\left\langle K_H^*\mathbb{D}\langle\vartheta(T),e^i\rangle, K_H^*\left(0,R_H^{-1}\left(\int_0^\cdot\sigma^*(X_u)e^i\d u\right)\right)\right\rangle_{L^2([0,T],\R^{d_1+d_2})}\cr
&=&\left\langle\mathbb{D}\langle\vartheta(T),e^i\rangle,\left(0,R_H^{-1}\left(\int_0^\cdot\sigma^*(X_u)e^i\d u\right)\right)\right\rangle_\mathcal {H}\cr
&=&\mathbb{D}_{\left(0,\int_0^\cdot\sigma^*(X_u)e^i\d u\right)}\langle\vartheta(T),e^i\rangle\cr
&=&\left\langle\left(\int_0^T(\sigma\sigma^*)(X_u)\d u\right)^{-1}\int_0^T\frac{T-u}{T}\nabla\sigma(X_u)v_1\sigma^*(X_u)e^i\d u,e^i\right\rangle.
\eeqlb
So, by (H2) it follows that the left hand of \eqref{Pf of The 3.1-14} is square integrable.
Combining this with \eqref{Pf of The 3.1-13}, by \cite[Proposition 1.3.3]{Nualart06a} we show that, for each $i=1,\cdots,d_2$,
\beqlb\label{Pf of The 3.1-15}
\langle\vartheta(T),e^i\rangle\left(0,K_H^{-1}\left(\int_0^\cdot\sigma^*(X_u)e^i\d u\right)(t)\right)\in\mathrm{Dom}\delta_{\bar{W}},
\eeqlb
and moreover,
\beqlb\label{Pf of The 3.1-16}
&&\delta_{\bar{W}}\left(\langle\vartheta(T),e^i\rangle\left(0,K_H^{-1}\left(\int_0^\cdot\sigma^*(X_u)e^i\d u\right)\right)\right)\cr
&=&\langle\vartheta(T),e^i\rangle\delta_{\bar{W}}\left(0,K_H^{-1}\left(\int_0^\cdot\sigma^*(X_u)e^i\d u\right)\right)\cr
&&-\left\langle\mathbb{D}_{\bar{W}}\langle\vartheta(T),e^i\rangle,\left(0,K_H^{-1}\left(\int_0^\cdot\sigma^*(X_u)e^i\d u\right)\right)\right\rangle_{\mathcal {H}_{\bar{W}}}\cr
&=&\langle\vartheta(T),e^i\rangle\int_0^T\left\langle K_H^{-1}\left(\int_0^\cdot\sigma^*(X_u)e^i\d u\right)(t),\d \tilde{W}_t\right\rangle\cr
&&-\left\langle\left(\int_0^T(\sigma\sigma^*)(X_u)\d u\right)^{-1}\int_0^T\frac{T-u}{T}\nabla\sigma(X_u)v_1\sigma^*(X_u)e^i\d u,e^i\right\rangle,
\eeqlb
where the last relation is due to \eqref{Pf of The 3.1-11} and \eqref{Pf of The 3.1-14}.\\
Consequently, we derive that $K_H^*h\in\mathrm{Dom}\delta_{\bar{W}}$ and by \eqref{Pf of The 3.1-6},\eqref{Pf of The 3.1-7} and \eqref{Pf of The 3.1-16},
\beqlb\label{Pf of The 3.1-17}
&&\delta_{\bar{W}}(K_H^*h)=\delta_{\bar{W}}(K_H^*{h_1},0)+\delta_{\bar{W}}(0,K_H^*{h_2})\cr
&=&\int_0^T\left\langle K_H^{-1}\left(\frac{\cdot}{T}\right)(t)v_1,\d W_t\right\rangle
+\sum\limits_{i=1}^{d_2}\langle\vartheta(T),e^i\rangle\int_0^T\left\langle K_H^{-1}\left(\int_0^\cdot\sigma^*(X_u)e^i\d u\right)(t),\d \tilde{W}_t\right\rangle\cr
&&-\sum\limits_{i=1}^{d_2}\left\langle\left(\int_0^T(\sigma\sigma^*)(X_u)\d u\right)^{-1}\int_0^T\frac{T-u}{T}\nabla\sigma(X_u)v_1\sigma^*(X_u)e^i\d u,e^i\right\rangle\cr
&=&\int_0^T\left\langle K_H^{-1}\left(\frac{\cdot}{T}\right)(t)v_1,\d W_t\right\rangle
+\sum\limits_{i=1}^{d_2}\langle\vartheta(T),e^i\rangle\left\langle\int_0^T\left(K_H^{-1}\left(\int_0^\cdot\sigma^*(X_u)\d u\right)\right)^*(t)\d\tilde{W}_t,e^i\right\rangle\cr
&&-\sum\limits_{i=1}^{d_2}\left\langle\left(\int_0^T(\sigma\sigma^*)(X_u)\d u\right)^{-1}\int_0^T\frac{T-u}{T}\nabla\sigma(X_u)v_1\sigma^*(X_u)\d u\cdot e^i,e^i\right\rangle\cr
&=&\int_0^T\left\langle K_H^{-1}\left(\frac{\cdot}{T}\right)(t)v_1,\d W_t\right\rangle
+\left\langle\vartheta(T),\int_0^T\left(K_H^{-1}\left(\int_0^\cdot\sigma^*(X_u)\d u\right)\right)^*(t)\d\tilde{W}_t\right\rangle\cr
&&-Tr\left(\left(\int_0^T(\sigma\sigma^*)(X_u)\d u\right)^{-1}\int_0^T\frac{T-u}{T}\nabla\sigma(X_u)v_1\sigma^*(X_u)\d u\right).
\eeqlb
Then it follows by Proposition \ref{Pro 2.2} that $h\in\mathrm{Dom}\delta$ and $\delta(h)=\delta_{\bar{W}}(K_H^*h)$.
The proof is finished.
\fin

\bremark\label{Rem 3.1}
(i) If $H=\frac{1}{2}$, the inverse operator $K_H^{-1}$ appeared in Theorem \ref{The 3.1} can be viewed as the derivative operator.
By simple calculus, we know that our result covers that of \cite[Theorem 1.1]{Wang12a}.

(ii) By \eqref{2.1}, it is not difficult to verify that the right side of \eqref{Th3.1-0'} equals to the following relation
\beqlb\label{Rem 3.1-1}
M_T
&=&c_1\int_0^T\left\langle\frac{t^{\frac{1}{2}-H}}{T}v_1,\d W_t\right\rangle
+c_2\left\langle\vartheta(T),\int_0^Tt^{\frac{1}{2}-H}\sigma(x+B^H_t)\d\tilde{W}_t\right\rangle\cr
&&+c_3\left\langle\vartheta(T),\int_0^Tt^{H-\frac{1}{2}}
\int_0^t\frac{\sigma(x+B^H_t)-\sigma(x+B^H_r)}{(t-r)^{\frac{1}{2}+H}}r^{\frac{1}{2}-H}\d r\d\tilde{W}_t\right\rangle\cr
&&-Tr\left(\left(\int_0^T(\sigma\sigma^*)\left(x+B^H_u\right)\d u\right)^{-1}\int_0^T\frac{T-u}{T}\nabla\sigma\left(x+B^H_u\right)v_1\sigma^*(x+B^H_u)\d u\right),\nonumber
\eeqlb
where $c_i,i=1,2,3$ are three constants.
\eremark

Let us go back to the above proof and note that, for \eqref{Pf of The 3.1-2},
the choice of $h_2$ is not unique.
Now, we begin with the assumption (H3):\\
$\sigma$ is differentiable with bounded derivative,
and $\sigma\sigma^*$ is invertible such that $(\sigma\sigma^*)^{-1}$ is bounded and H\"{o}lder continuous of order $\gamma\in(1-1/(2H),1]$:
\beqnn
\|(\sigma\sigma^*)^{-1}(z_1)-(\sigma\sigma^*)^{-1}(z_2)\|\leq K|z_1-z_2|^\gamma, \ \ \forall z_1,z_2\in\R^{d_1},
\eeqnn
where $K$ is a positive constant.

Set
\begin{equation}\label{Pf of The 3.1-18}
h_2(t)
=R_H^{-1}\left(\int_0^\cdot\left(\sigma^*(\sigma\sigma^*)^{-1}\right)(X_u)\d u\cdot\frac{1}{T}\left(v_2+\int_0^T\frac{T-u}{T}\nabla\sigma(X_u)v_1\d \tilde{B}^H_u\right)\right)(t).
\end{equation}
It is not hard to verify that $h_1$ defined as \eqref{Pf of The 3.1-1} before and the above $h_2$ also satisfy \eqref{3.4'} and \eqref{3.4}.
In the spirit of the proof of Theorem \ref{The 3.1}, we have the following derivative formula for the equation \eqref{3.1}.
\btheorem\label{The 3.2}
Assume (H3) and let $v=(v_1,v_2)\in\R^{d_1+d_2}$.
Then we have
\beqlb\label{Th3.2-1}
\nabla_v P_T f(x,y)=\E^{x,y}\left[f(X_T,Y_T)\tilde{M}_T\right],\ \ f\in C_b^1(\R^{d_1+d_2}),
\eeqlb
where
\beqlb\label{Th3.2-2}
\tilde{M}_T
&=&\int_0^T\left\langle K_H^{-1}\left(\frac{\cdot}{T}\right)(t)v_1,\d W_t\right\rangle\cr
&&+\left\langle\tilde{\vartheta}(T),\int_0^T\left(K_H^{-1}\left(\int_0^\cdot\left(\sigma^*(\sigma\sigma^*)^{-1}\right)\left(x+B^H_u\right)\d u\right)\right)^*(t)\d\tilde{W}_t\right\rangle\cr
&&-Tr\left(\int_0^T\frac{T-u}{T^2}\nabla\sigma\left(x+B^H_u\right)v_1\left(\sigma^*(\sigma\sigma^*)^{-1}\right)\left(x+B^H_u\right)\d u\right)
\eeqlb
and
\beqlb\label{Th3.2-3}
\tilde{\vartheta}(T)=\frac{1}{T}\left(v_2+\int_0^T\frac{T-u}{T}\nabla\sigma\left(x+B^H_u\right)v_1\d \tilde{B}^H_u\right).
\eeqlb
\etheorem

\bremark\label{Rem 3.2}
Theorem \ref{The 3.1} and Theorem \ref{The 3.2} are established under the conditions (H2) and (H3), respectively.
Observe that these two assumptions are not compatible.
For instance, setting
\begin{equation}\nonumber
\zeta(t)=\left(
\begin{array}{ccc}
1 & t\\
t & t^2\\
\end{array}
\right),
\end{equation}
it is easy to see that $\int_0^T\zeta(t)\d t$ is invertible, yet $\zeta(t)$ is degenerate.
On the other hand, taking

\begin{equation}\nonumber
\eta(t)=\left(
\begin{array}{ccc}
\sin t & \cos t\\
\cos t & -\sin t\\
\end{array}
\right),
\end{equation}
it is not hard to verify that $\eta(t)$ is invertible, while $\int_0^T\eta(t)\d t$ is degenerate.
\eremark

The next result states that with the help of Theorem \ref{The 3.2}, an explicit gradient estimate can be derived.
\bcorollary\label{Cor 3.1}
Assume (H3).
Then for each $p>1/\left(\frac{3}{2}-H\right)$ there exists a constant $C(p,H)>0$ such that
\beqlb\label{Cor 3.1-1}\nonumber
|\nabla_v P_T f(x,y)|&\leq&C(p,H)(P_T|f|^p)^\frac{1}{p}(x,y)\cr
&&\times\left[|v_1|\left(\frac{1}{T^{H}}+1+T^{\gamma(H-\tilde{\epsilon})}+T^{H-\tilde{\epsilon}}\right)
+|v_2|\left(\frac{1}{T^{H}}+\frac{1}{T^{H-\gamma(H-\tilde{\epsilon})}}\right)\right]
\eeqlb
holds for $v=(v_1,v_2)\in\R^{d_1+d_2}$ and $(x,y)\in\R^{d_1+d_2}$,
where $\tilde{\epsilon}$ is fixed with $\tilde{\epsilon}<H-\frac{1}{\gamma}\left(H-\frac{1}{2}\right)$.
\ecorollary

Before proving Corollary \ref{Cor 3.1}, we will first present a technical lemma used to calculate the moment of stochastic integral
w.r.t. a Brownian motion (see \cite[Lemma 2.2]{Wang12a}).

\blemma\label{Lem 3.2}
Let $\varrho(t)$ be a predictable process on $\R^d$ with $\E\int_0^T|\varrho(t)|^q\d t<\infty$ for some $q\geq 2$.
Then there holds
\beqnn
\E\left|\int_0^T\langle\varrho(t),\d\check{W}_t\rangle\right|^q\leq CT^{\frac{q-2}{2}}\int_0^T\E|\varrho(t)|^q\d t,
\eeqnn
where $\check{W}_t$ is a Brownian motion on $\R^d$.
\elemma

\emph{Proof of Corollary \ref{Cor 3.1}.}
We first note that by the H\"{o}lder inequality, it suffices to prove for $p\in(1,2]$.
According to Theorem \ref{The 3.2} and the H\"{o}lder inequality again, we have
\beqlb\label{Pf of Cor 3.1-0}
|\nabla_v P_T f(x,y)|=\left|\E^{x,y}\left[f(X_T,Y_T)\tilde{M}_T\right]\right|\leq (P_T|f|^p)^\frac{1}{p}(x,y)\cdot\left(\E|\tilde{M}_T|^q\right)^\frac{1}{q},
\eeqlb
where $q=\frac{p}{p-1}$.\\
Now, we shall bound the term $\E|\tilde{M}_T|^q$.
For the convenience of the notations, let
$$\psi(z)=(\sigma^*(\sigma\sigma^*)^{-1})(z),\ \ z\in\R^{d_1}.$$
By \eqref{2.1}, we first obtain
\beqlb\label{Pf of Cor 3.1-1}
&&K_H^{-1}\left(\frac{\cdot}{T}\right)(t)
=t^{H-\frac{1}{2}}D_{0+}^{H-\frac{1}{2}}\left(\cdot^{\frac{1}{2}-H}\frac{1}{T}\right)(t)\cr
&=&\frac{1}{\Gamma(\frac{3}{2}-H)}\left[\frac{t^{\frac{1}{2}-H}}{T}+
+\left(H-\frac{1}{2}\right)\frac{t^{H-\frac{1}{2}}}{T}\int_0^t\frac{t^{\frac{1}{2}-H}-r^{\frac{1}{2}-H}}{(t-r)^{\frac{1}{2}+H}}\d r\right]\cr
&=&\frac{1-C_0(H-\frac{1}{2})}{\Gamma(\frac{3}{2}-H)}\frac{t^{\frac{1}{2}-H}}{T}
\eeqlb
and
\beqlb\label{Pf of Cor 3.1-2}
&&K_H^{-1}\left(\int_0^\cdot\psi\left(x+B^H_u\right)\d u\right)(t)
=t^{H-\frac{1}{2}}D_{0+}^{H-\frac{1}{2}}\left(\cdot^{\frac{1}{2}-H}\psi(x+B^H_\cdot)\right)(t)\cr
&=&\frac{1}{\Gamma(\frac{3}{2}-H)}\Bigg[t^{\frac{1}{2}-H}\psi(x+B^H_t)
+\left(H-\frac{1}{2}\right)\psi(x+B^H_t)t^{H-\frac{1}{2}}\int_0^t\frac{t^{\frac{1}{2}-H}-r^{\frac{1}{2}-H}}{(t-r)^{\frac{1}{2}+H}}\d r\cr
&&+\left(H-\frac{1}{2}\right)t^{H-\frac{1}{2}}
\int_0^t\frac{\psi(x+B^H_t)-\psi(x+B^H_r)}{(t-r)^{\frac{1}{2}+H}}r^{\frac{1}{2}-H}\d r\Bigg]\cr
&=&\frac{1-C_0(H-\frac{1}{2})}{\Gamma(\frac{3}{2}-H)}t^{\frac{1}{2}-H}\psi(x+B^H_t)\cr
&&+\frac{H-\frac{1}{2}}{\Gamma(\frac{3}{2}-H)}t^{H-\frac{1}{2}}
\int_0^t\frac{\psi(x+B^H_t)-\psi(x+B^H_r)}{(t-r)^{\frac{1}{2}+H}}r^{\frac{1}{2}-H}\d r.
\eeqlb
Recall that the last equalities of \eqref{Pf of Cor 3.1-1} and \eqref{Pf of Cor 3.1-2} are due to \eqref{Pf of The 3.1-9'}.
Combining \eqref{Pf of Cor 3.1-1} with Lemma \ref{Lem 3.2} yields
\beqlb\label{Pf of Cor 3.1-3}
\E\left|\int_0^T\left\langle K_H^{-1}\left(\frac{\cdot}{T}\right)(t)v_1,\d W_t\right\rangle\right|^q
=C\E\left|\int_0^T\left\langle\frac{t^{\frac{1}{2}-H}}{T}v_1,\d W_t\right\rangle\right|^q
\leq C|v_1|^q\frac{1}{T^{Hq}}.
\eeqlb

Next, we are to estimate the $q$th moment of the second term of the right side of \eqref{Th3.2-2}.
By \eqref{Pf of Cor 3.1-2} and (H3), we first get, for each $1<\alpha<1/\left(H-\frac{1}{2}\right)$,
\beqlb\label{Pf of Cor 3.1-4'}\nonumber
&&\int_0^T\left|K_H^{-1}\left(\int_0^\cdot
\psi\left(x+B^H_u\right)\d u\right)(t)\right|^\alpha\d t\cr
&\leq&C\Bigg(\int_0^Tt^{\left(\frac{1}{2}-H\right)\alpha}\left|\psi(x+B^H_t)\right|^\alpha\d t
+\int_0^T
\left|t^{H-\frac{1}{2}}\int_0^t\frac{\psi(x+B^H_t)-\psi(x+B^H_r)}{(t-r)^{\frac{1}{2}+H}}r^{\frac{1}{2}-H}\d r\right|^\alpha\d t
\Bigg)\cr
&\leq&C\left(T^{\left(\frac{1}{2}-H\right)\alpha+1}
+\int_0^T
\left(t^{H-\frac{1}{2}}\int_0^t\frac{|B^H_t-B^H_r|^\gamma+|B^H_t-B^H_r|}{(t-r)^{\frac{1}{2}+H}}r^{\frac{1}{2}-H}\d r\right)^\alpha\d t\right)\cr
&\leq&C\left(T^{\left(\frac{1}{2}-H\right)\alpha+1}
+\int_0^T
t^{\left(\gamma(H-\tilde{\epsilon})+\frac{1}{2}-H\right)\alpha}\d t\cdot\|B^H\|_{H-\tilde{\epsilon}}^{\gamma \alpha}
+\int_0^T
t^{\left(\frac{1}{2}-\tilde{\epsilon}\right)\alpha}\d t\cdot\|B^H\|_{H-\tilde{\epsilon}}^\alpha
\right)\cr
&=&C\left(T^{\left(\frac{1}{2}-H\right)\alpha+1}+T^{\left(\gamma(H-\tilde{\epsilon})+\frac{1}{2}-H\right)\alpha+1}\|B^H\|_{H-\tilde{\epsilon}}^{\gamma \alpha}
+T^{\left(\frac{1}{2}-\tilde{\epsilon}\right)\alpha+1}\cdot\|B^H\|_{H-\tilde{\epsilon}}^\alpha
\right),
\eeqlb
where $\tilde{\epsilon}$ is fixed with $\tilde{\epsilon}<H-\frac{1}{\gamma}\left(H-\frac{1}{2}\right)$.\\
Then, using the fact that $x+B^H_t$ is measurable w.r.t. $\mathcal{F}_T$ while $\tilde{W}$ is independent of $\mathcal{F}_T$,
Lemma \ref{Lem 3.2} and \cite[Theorem 1.1]{Memin&Mishura&Valkeila01} along with (H3), we derive that,
for each $1<q<1/\left(H-\frac{1}{2}\right) \left(\mathrm{or}\ p>1/\left(\frac{3}{2}-H\right)\right)$,
\beqlb\label{Pf of Cor 3.1-4}\nonumber
&&\E\left(\left|\left\langle\frac{v_2}{T},\int_0^T\left(K_H^{-1}\left(\int_0^\cdot
\psi\left(x+B^H_u\right)\d u\right)\right)^*(t)\d\tilde{W}_t\right\rangle\right|^q|\mathcal{F}_T\right)\cr
&\leq&C\frac{|v_2|^q}{T^{\frac{q}{2}+1}}\int_0^T\left|K_H^{-1}\left(\int_0^\cdot
\psi\left(x+B^H_u\right)\d u\right)(t)\right|^q\d t\cr
&\leq&C\frac{|v_2|^q}{T^{\frac{q}{2}+1}}\left(T^{\left(\frac{1}{2}-H\right)q+1}
+T^{\left(\gamma(H-\tilde{\epsilon})+\frac{1}{2}-H\right)q+1}\|B^H\|_{H-\tilde{\epsilon}}^{\gamma q}+
T^{\left(\frac{1}{2}-\tilde{\epsilon}\right)q+1}\cdot\|B^H\|_{H-\tilde{\epsilon}}^q
\right)\cr
&=&C|v_2|^q\left(\frac{1}{T^{Hq}}+\frac{\|B^H\|_{H-\tilde{\epsilon}}^{\gamma q}}{T^{\left(H-\gamma(H-\tilde{\epsilon})\right)q}}
+\frac{\|B^H\|_{H-\tilde{\epsilon}}^{q}}{T^{\tilde{\epsilon}q}}
\right)
\eeqlb
and
\beqlb\label{Pf of Cor 3.1-5}\nonumber
&&\E\left(\left|\left\langle\int_0^T\frac{T-u}{T^2}\nabla\sigma\left(x+B^H_u\right)v_1\d \tilde{B}^H_u,
\int_0^T\left(K_H^{-1}\left(\int_0^\cdot\psi\left(x+B^H_u\right)\d u\right)\right)^*(t)\d\tilde{W}_t\right\rangle\right|^q|\mathcal{F}_T\right)\cr
&\leq&
\left[\E\left(\left|\int_0^T\frac{T-u}{T^2}\nabla\sigma\left(x+B^H_u\right)v_1\d \tilde{B}^H_u\right|^{2q}|\mathcal{F}_T\right)\right]^\frac{1}{2}\cr
&&\times\left[\E\left(\left|\int_0^T\left(K_H^{-1}\left(\int_0^\cdot\psi\left(x+B^H_u\right)\d u\right)\right)^*(t)\d\tilde{W}_t\right|^{2q}|\mathcal{F}_T\right)\right]^\frac{1}{2}\cr
&\leq&C\left(\int_0^T\left|\frac{T-u}{T^2}\nabla\sigma\left(x+B^H_u\right)v_1\right|^{\frac{1}{H}}\d u\right)^{Hq}
\cdot\left(T^{q-1}\int_0^T\left|K_H^{-1}\left(\int_0^\cdot\psi\left(x+B^H_u\right)\d u\right)^*(t)\right|^{2q}\d t\right)^\frac{1}{2}\cr
&\leq&C|v_1|^q\frac{1}{T^{(1-H)q}}T^{\frac{q-1}{2}}
\left(T^{\left(\frac{1}{2}-H\right)q+\frac{1}{2}}+T^{\left(\gamma(H-\tilde{\epsilon})+\frac{1}{2}-H\right)q+\frac{1}{2}}\|B^H\|_{H-\tilde{\epsilon}}^{\gamma q}
+T^{\left(\frac{1}{2}-\tilde{\epsilon}\right)q+\frac{1}{2}}\cdot\|B^H\|_{H-\tilde{\epsilon}}^q\right)\cr
&=&C|v_1|^q\left(1+T^{\gamma(H-\tilde{\epsilon})q}\|B^H\|_{H-\tilde{\epsilon}}^{\gamma q}
+T^{(H-\tilde{\epsilon})q}\|B^H\|_{H-\tilde{\epsilon}}^q
\right).
\eeqlb
Consequently, these, together with \eqref{Th3.2-3}, lead to
\beqlb\label{Pf of Cor 3.1-6}
&&\E\left|\left\langle\tilde{\vartheta}(T),\int_0^T\left(K_H^{-1}\left(\int_0^\cdot
\left(\sigma^*(\sigma\sigma^*)^{-1}\right)\left(x+B^H_u\right)\d u\right)\right)^*(t)\d\tilde{W}_t\right\rangle\right|^q\cr
&\leq&
C\E\left[\E\left(\left|\left\langle\frac{v_2}{T},\int_0^T\left(K_H^{-1}\left(\int_0^\cdot
\psi\left(x+B^H_u\right)\d u\right)\right)^*(t)\d\tilde{W}_t\right\rangle\right|^q|\mathcal{F}_T\right)\right]\cr
&&+
C\E\left[\E\left(\left|\left\langle\int_0^T\frac{T-u}{T^2}\nabla\sigma\left(x+B^H_u\right)v_1\d \tilde{B}^H_u,
\int_0^T\left(K_H^{-1}\left(\int_0^\cdot\psi\left(x+B^H_u\right)\d u\right)\right)^*(t)\d\tilde{W}_t\right\rangle\right|^q|\mathcal{F}_T\right)\right]\cr
&\leq&C|v_2|^q\left(\frac{1}{T^{Hq}}+\frac{1}{T^{\left(H-\gamma(H-\tilde{\epsilon})\right)q}}\right)
+C|v_1|^q\left(1+T^{\gamma(H-\tilde{\epsilon})q}+T^{(H-\tilde{\epsilon})q}\right).
\eeqlb
As for the third term of the right side of \eqref{Th3.2-2}, by (H3) it is easy to see that
\beqlb\label{Pf of Cor 3.1-7}
&&\E\left|Tr\left(\int_0^T\frac{T-u}{T^2}\nabla\sigma\left(x+B^H_u\right)v_1\left(\sigma^*(\sigma\sigma^*)^{-1}\right)\left(x+B^H_u\right)\d u\right)\right|^q\cr
&\leq&C\frac{|v_1|^q}{T^q}\E\left(\int_0^T\left|\nabla\sigma\left(x+B^H_u\right)\right|\d u\right)^q\cr
&\leq&C|v_1|^q.
\eeqlb
So, by \eqref{Pf of Cor 3.1-3}-\eqref{Pf of Cor 3.1-7} we get
\beqlb\label{Pf of Cor 3.1-8}\nonumber
\E|\tilde{M}_T|^q&\leq&C\E\left|\int_0^T\left\langle K_H^{-1}\left(\frac{\cdot}{T}\right)(t)v_1,\d W_t\right\rangle\right|^q\cr
&&+C\E\left|\left\langle\tilde{\vartheta}(T),\int_0^T\left(K_H^{-1}\left(\int_0^\cdot
\left(\sigma^*(\sigma\sigma^*)^{-1}\right)\left(x+B^H_u\right)\d u\right)\right)^*(t)\d\tilde{W}_t\right\rangle\right|^q\cr
&&+C\E\left|Tr\left(\int_0^T\frac{T-u}{T^2}\nabla\sigma\left(x+B^H_u\right)v_1\left(\sigma^*(\sigma\sigma^*)^{-1}\right)\left(x+B^H_u\right)\d u\right)\right|^q\cr
&\leq&C|v_1|^q\frac{1}{T^{Hq}}
+C|v_2|^q\left(\frac{1}{T^{Hq}}+\frac{1}{T^{\left(H-\gamma(H-\epsilon)\right)q}}\right)\cr
&&+C|v_1|^q\left(1+T^{\gamma(H-\tilde{\epsilon})q}+T^{(H-\tilde{\epsilon})q}\right)
+C|v_1|^q\cr
&=&C|v_1|^q\left(\frac{1}{T^{Hq}}+1+T^{\gamma(H-\tilde{\epsilon})q}+T^{(H-\tilde{\epsilon})q}\right)
+C|v_2|^q\left(\frac{1}{T^{Hq}}+\frac{1}{T^{\left(H-\gamma(H-\epsilon)\right)q}}\right).\cr
\eeqlb
Then, the desired assertion follows from \eqref{Pf of Cor 3.1-0}.
\fin

\bremark\label{Rem 3.3}
If we use the approach as in the proof of Corollary \ref{Cor 3.1} to obtain gradient estimate from Theorem \ref{The 3.1},
there will appear the following term
\beqnn
\left\|\left(\int_0^T(\sigma\sigma^*)\left(x+B^H_u\right)\d u\right)^{-1}\right\|
\eeqnn
which seems to be very difficult to deal with.
So, compared with Theorem \ref{The 3.1},
the advantage of Theorem \ref{The 3.2} is that an explicit gradient estimate is obtained with simpler proofs,
while the drawback is the restriction to non-degenerate noise situation.
\eremark

\section{An extension to a more general model}

\setcounter{equation}{0}

This section is devoted to extend Theorem \ref{The 3.1} to a more general model.
That is,
\begin{equation}\label{4.1}
\begin{cases}
 \textnormal\d X_t=b_1(X_t)\d t+\sigma_1\d B_t^H,\ \ \ \ \ \ \ \ \ \ \ \ \ \ \ \ \ X_0=x\in\R^{d_1},\\
 \textnormal\d Y_t=b_2(X_t)\d t+\sigma_2(X_t)\d\tilde{B}_t^H, \ \ \ \ \ \ \ \ \ \ \ \ Y_0=y\in\R^{d_2},
\end{cases}
\end{equation}
where $b_1:\R^{d_1}\rightarrow\R^{d_1}, b_2:\R^{d_1}\rightarrow\R^{d_2}, \sigma_1\in\R^{d_1}\times\R^{d_1}$ is invertible,
$\sigma_2:\R^{d_1}\rightarrow\R^{d_2}\times\R^l$
and $(B_t^H,\tilde{B}_t^H)$ is a fractional Brownian motion on $\R^{d_1+l}$.
In this setting, one can allow $\sigma_2$ to be degenerate.
Let $P_t$ be the associated operators.
To establish the derivative formula for $P_t$,
we shall make use of the following assumption.\\
(H4)
$b_i,i=1,2$ and $\sigma_2$ are differentiable with bounded derivatives, and $\int_0^T(\sigma_2\sigma_2^*)\left(X_u\right)\d u$
is invertible such that
\beqnn
\E\left\|\left(\int_0^T(\sigma_2\sigma_2^*)\left(X_u\right)\d u\right)^{-1}\right\|^{2+\epsilon_0}<\infty,
\eeqnn
where $\epsilon_0>0$ is a fixed constant.

Obviously, by (H4), it follows from \cite[Theorem 2.1]{Nualart&Rascanu02a} that \eqref{4.1} has a unique solution.
The main result of this part is the following.
\btheorem\label{The 4.1}
Assume (H4) and let $v=(v_1,v_2)\in\R^{d_1+d_2}$.
Then
\beqlb\label{Th4.1-0}
\nabla_v P_T f(x,y)=\E^{x,y}\left[f(X_T,Y_T)N_T\right],\ \ f\in C_b^1(\R^{d_1+d_2}),
\eeqlb
where
\beqlb\label{Th4.1-1}
N_T
&=&\int_0^T\left\langle \sigma_1^{-1}K_H^{-1}\left(\frac{T-\cdot}{T}\nabla b_1(X_\cdot)+\frac{1}{T}\right)(t)v_1,\d W_t\right\rangle\cr
&&+\left\langle\chi(T),\int_0^T\left(K_H^{-1}\left(\int_0^\cdot\sigma_2^*(X_u)\d u\right)\right)^*(t)\d\tilde{W}_t\right\rangle\cr
&&-Tr\left(\left(\int_0^T(\sigma_2\sigma_2^*)(X_u)\d u\right)^{-1}\int_0^T\frac{T-u}{T}\nabla\sigma_2(X_u)v_1\sigma_2^*(X_u)\d u\right)
\eeqlb
and
\beqlb\label{Th4.1-2}
\chi(T)=\left(\int_0^T(\sigma_2\sigma_2^*)(X_u)\d u\right)^{-1}
\left(v_2+\int_0^T\frac{T-u}{T}\nabla b_2(X_u)v_1\d u
+\int_0^T\frac{T-u}{T}\nabla\sigma_2(X_u)v_1\d \tilde{B}^H_u\right).\cr
\eeqlb
\etheorem
\emph{Proof.}
Let $h=(h_1,h_2)$ be as follows: for $t\in[0,T]$,
\beqlb\label{Pf of The 4.1-1}
h_1(t)=R_H^{-1}\left(\int_0^\cdot\sigma_1^{-1}\left(\frac{T-u}{T}\nabla b_1(X_u)+\frac{1}{T}\right)v_1\d u\right)(t)
\eeqlb
and
\beqlb\label{Pf of The 4.1-2}
h_2(t)&=&R_H^{-1}\Bigg(\int_0^\cdot\sigma_2^*(X_u)\d u\cdot\left(\int_0^T(\sigma_2\sigma_2^*)(X_u)\d u\right)^{-1}\cr
&&\ \ \ \ \ \ \ \ \cdot\left(v_2+\int_0^T\frac{T-u}{T}\nabla b_2(X_u)v_1\d u
+\int_0^T\frac{T-u}{T}\nabla\sigma_2(X_u)v_1\d \tilde{B}^H_u\right)\Bigg)(t).\nonumber\\
\eeqlb
Following the argument of Theorem \ref{The 3.1},
by (H4) we show that the above $h\in{\rm Dom}\delta$ and then $\delta(h)=N_T$.
Hence, if one has $(\mathbb{D}_{R_Hh}X_T,\mathbb{D}_{R_Hh}Y_T)=(\nabla_vX_T,\nabla_vY_T)$,
then by \eqref{Pf of Lem 3.2-6} we get the desired assertion.

Indeed, by Lemma \ref{Lem 3.1} and the definition of $h_1$ defined as \eqref{Pf of The 4.1-1} we have
\beqlb\label{Pf of Tem 4.1-3}\nonumber
\mathbb{D}_{R_Hh}X_t&=&\int_0^t\nabla b_1(X_u)\mathbb{D}_{R_Hh}X_u\d u+\int_0^t\sigma_1(R_H h_1)'(t)\d u\cr
&=&\int_0^t\nabla b_1(X_u)\mathbb{D}_{R_Hh}X_u\d u+\int_0^t\left(\frac{T-u}{T}\nabla b_1(X_u)+\frac{1}{T}\right)v_1\d u.
\eeqlb
This leads to
\beqlb\label{Pf of Tem 4.1-4}\nonumber
\mathbb{D}_{R_Hh}X_t+\frac{T-t}{T}v_1
=v_1+\int_0^t\nabla b_1(X_u)\left(\mathbb{D}_{R_Hh}X_u+\frac{T-u}{T}v_1\right)\d u.
\eeqlb
Observe that $\nabla_v X_t$ satisfies the same equation, i.e.
\beqlb\label{Pf of Tem 4.1-5}\nonumber
\nabla_v X_t
=v_1+\int_0^t\nabla b_1(X_u)\nabla_v X_u\d u.
\eeqlb
Then by the uniqueness of the ODE we obtain
\beqlb\label{Pf of Tem 4.1-6}
\mathbb{D}_{R_Hh}X_t+\frac{T-t}{T}v_1=\nabla_v X_t,\ \ t\in[0,T].
\eeqlb
In particular, there holds $\mathbb{D}_{R_Hh}X_T=\nabla_v X_T$.\\
As for $\mathbb{D}_{R_Hh}Y_t$ and $\nabla_v Y_t$, we have
\beqlb\label{Pf of Tem 4.1-7}\nonumber
\mathbb{D}_{R_Hh}Y_t=\int_0^t\nabla b_2(X_u)\mathbb{D}_{R_Hh}X_u\d u
+\int_0^t\nabla\sigma_2(X_u)\mathbb{D}_{R_Hh}X_u\d\tilde{B}^H_u
+\int_0^t\sigma_2(X_u)(R_H h_2)'(u)\d u
\eeqlb
and
\beqlb\label{Pf of Tem 4.1-8}\nonumber
\nabla_v Y_t=v_2+
\int_0^t\nabla b_2(X_u)\nabla_v X_u\d u
+\int_0^t\nabla\sigma_2(X_u)\nabla_v X_u\d\tilde{B}^H_u.
\eeqlb
Consequently, by \eqref{Pf of Tem 4.1-6} we get
 \beqlb\label{Pf of Tem 4.1-9}\nonumber
\mathbb{D}_{R_Hh}Y_T-\nabla_v Y_T
&=&-v_2+\int_0^T\nabla b_2(X_u)\left(\mathbb{D}_{R_Hh}X_u-\nabla_v X_u\right)\d u\cr
&&+\int_0^T\nabla\sigma_2(X_u)\left(\mathbb{D}_{R_Hh}X_u-\nabla_v X_u\right)\d\tilde{B}^H_u
+\int_0^T\sigma_2(X_u)(R_H h_2)'(u)\d u\cr
&=&-v_2-\int_0^T\frac{T-u}{T}\nabla b_2(X_u)v_1\d u\cr
&&-\int_0^T\frac{T-u}{T}\nabla\sigma_2(X_u)v_1\d\tilde{B}^H_u+\int_0^T\sigma_2(X_u)(R_H h_2)'(u)\d u.
\eeqlb
Therefore, it follows from \eqref{Pf of The 4.1-2} that $\mathbb{D}_{R_Hh}Y_T=\nabla_v Y_T$.
This completes the proof.
\fin

\end{document}